\documentclass[preprint,11pt]{article}

\usepackage{lineno}
\usepackage{amssymb}
\usepackage{amsmath}
\usepackage{amsthm}
\usepackage{caption}
\usepackage{subcaption}
\usepackage{algorithm}
\usepackage{algpseudocode}
\usepackage{xcolor}
\usepackage{enumitem}
\usepackage{etoolbox}
\usepackage{comment}

\usepackage[hidelinks,
colorlinks,linkcolor=blue,citecolor=blue,final]{hyperref}

\let\oldenumerate\enumerate
\renewcommand{\enumerate}{
	\oldenumerate
	\setlength{\itemsep}{1.5pt}
	\setlength{\parskip}{0pt}
	\setlength{\parsep}{0pt}
}
\newtheorem{theorem}{Theorem}[section]

\newtheorem{conj}{Conjecture}
\newtheorem{quest}{Question}
\newtheorem{obs}{Observation}[section]
\usepackage{graphicx}
  \baselineskip 8 mm

\begin{document}
	%\begin{frontmatter}
	
	\title{Locating-dominating partitions for some classes of graphs}
	
	\author{$^1$Florent Foucaud, $^{2}$Paras Vinubhai Maniya, $^3$Kaustav Paul, and $^2$Dinabandhu Pradhan\thanks{Corresponding author.}  \\ \\
		$^{1}$Université Clermont Auvergne, CNRS, Clermont Auvergne INP,\\  Mines Saint-Étienne, LIMOS, 63000 Clermont-Ferrand, France\\ \\
		\small \tt Email: florent.foucaud@uca.fr
		\\ \\
		$^{2}$Department of Mathematics \& Computing\\
		Indian Institute of Technology (ISM) \\
		Dhanbad, India \\
		\small \tt Email: maniyaparas9999@gmail.com \\
		\small \tt Email: dina@iitism.ac.in\\
		\\\\
		$^{3}$Department of Mathematics, \\
		Indian Institute of Technology Ropar \\
		Rupnagar, India \\
		\small \tt Email: kaustav.20maz0010@iitrpr.ac.in}
	
	\date{}
	\maketitle

%\author[inst2]{}

\begin{abstract}
A dominating set of a graph $G$ is a set $D \subseteq V(G)$ such that every vertex in $V(G) \setminus D$ is adjacent to at least one vertex in $D$. A set $L\subseteq V(G)$ is a locating set of $G$ if every two vertices in $V(G) \setminus L$ have pairwise distinct open neighborhoods in $L$. A set $D\subseteq V(G)$ is a locating-dominating set of $G$ if $D$ is a dominating set and a locating set of $G$. The location-domination number of $G$, denoted by $\gamma_{LD}(G)$, is the minimum cardinality among all locating-dominating sets of $G$. A well-known conjecture in the study of locating-dominating sets is that if $G$ is an isolate-free and twin-free graph of order $n$, then $\gamma_{LD}(G)\le \frac{n}{2}$. Recently, Bousquet et al. [Discrete Math. 348 (2025), 114297] proved that if $G$ is an isolate-free and twin-free graph of order $n$, then $\gamma_{LD}(G)\le \lceil\frac{5n}{8}\rceil$ and posed the question whether the vertex set of such a graph can be partitioned into two locating sets. We answer this question affirmatively for twin-free distance-hereditary graphs, maximal outerplanar graphs, split graphs, and co-bipartite graphs.  In fact, we prove a stronger result: for any graph $G$ without isolated vertices and twin vertices, if $G$ is a distance-hereditary graph or a maximal outerplanar graph or a split graph or a co-bipartite graph, then the vertex set of $G$ can be partitioned into two locating-dominating sets. Consequently, this also confirms the original conjecture for these graph classes.

\end{abstract}



{\small \textbf{Keywords:} Domination number; location-domination number;  LD-partition;  distance-hereditary graphs;  maximal outerplanar graphs;    split graphs;  co-bipartite graphs.} \\
\indent {\small \textbf{AMS subject classification: 05C69}}

%\end{frontmatter}

\section{Introduction}
\label{sec:intro}
All the graphs considered in this paper are finite, simple, and undirected. For a graph $G$, we use $V(G)$ and $E(G)$ to denote the vertex set and the edge set of $G$, respectively. Two vertices $u$ and $v$ of $G$ are \emph{adjacent} if $uv\in E(G)$. The $\emph{neighbors}$ of $v$ in $G$ are the vertices adjacent to $v$ in $G$. The \emph{open neighborhood} $N_G(v)$ of a vertex $v$ in $G$ is the set of neighbors of $v$, while the \emph{closed neighborhood} of $v$ is the set $N_G[v] = \{v\} \cup N_G(v)$. The \emph{degree} of a vertex $v$ in $G$ is the number of vertices adjacent to $v$ in $G$, and is denoted by $\deg_G(v)$. An isolated vertex in a graph $G$ is a vertex of degree $0$. A graph without any isolated vertex is called an \emph{isolate-free} graph. A vertex of degree $1$ in $G$ is a \emph{leaf} of $G$. The \emph{maximum degree} of $G$ is the value $\max\{\deg_G(v) \colon v\in V(G)\}$. For vertices $u$ and $v$, $u$ and $v$ are called false (respectively, true) twins in $G$ if $N_G(u)=N_G(v)$ (respectively, $N_G[u]=N_G[v]$). Further, $u$ and $v$ are \emph{twins} in $G$ if they are false twins or true twins in $G$. A graph is \emph{twin-free} if it does not contain twins. In a rooted tree, one vertex is designated as the root. Consider a tree $T$ with the vertex $r$ as the root. For each vertex $v\neq r$ of $T$, the \emph{parent} of $v$ in $T$ is the neighbor of $v$ on the unique path from the root $r$ to $v$ in $T$. A \emph{child} of $v$ in $T$ is any of its neighbors other than its parent.

In this paper, we study distance-hereditary graphs, maximal outerplanar graphs, split graphs, and co-bipartite graphs. A \emph{distance-hereditary graph} is a graph in which the distance between any two vertices in any connected induced subgraph is the same as in the original graph. Their structure allows them to be built up recursively, which makes them useful for studying certain domination parameters (see ~\cite{Banerjee23,Liedloff08,Paul25}). Similarly, maximal outerplanar graphs, abbreviated as mops form a fundamental subclass of planar graphs. A graph is a \emph{mop} if it can be embedded in the plane such that all vertices lie on the boundary of its outer face (unbounded face) and all interior faces are triangles. Their well-understood structure allows for detailed combinatorial analysis, and they have been extensively studied in the context of domination parameters (see~\cite{Alvarado18, Araki18, Chvátal75,Claverol21,Dorfling6}). The set $X\subseteq V(G)$ is called a \emph{clique} (\emph{independent set}) of $G$ if every pair of vertices of $X$ are adjacent (nonadjacent) in $G$. A graph is a \emph{split graph} if its  vertex set can be partitioned into an independent set and a clique. A graph is a \emph{co-bipartite graph} if its vertex set can be partitioned into two cliques.
		
A dominating set $D$ of $G$ is a \emph{locating-dominating set}, abbreviated as LD-set, of $G$ if all vertices not in $D$ have pairwise distinct open neighborhoods in $D$. In other words, for every pair of vertices $u,v\in V(G)\setminus D$, we have $N_G(u)\cap D\neq N_G(v)\cap D$. The \emph{location-domination number} of $G$, denoted by $\gamma_{LD}(G)$, is the minimum cardinality among all LD-sets of $G$. Slater \cite{Slater88} in 1988 introduced this variant of domination, namely location-domination. Since its birth, location-domination remained an active area of research (see \cite{Bousquet25,Foucaud16-1,Foucaud17, Foucaud16-2,Garijo14}). This is due to its relevance in network science and theoretical computer science. For a comprehensive overview of locating-dominating sets in graphs, we recommend the book chapter~\cite{Lobstein20}.

\subsection{Motivation}

Research on locating-dominating sets has been significantly influenced by a conjecture made by Garijo et al.~\cite{Garijo14} in 2014. Foucaud and Henning~\cite{Foucaud16-1} later proposed a reformulation of this conjecture. The conjecture is stated below.

\begin{conj}[\cite{Foucaud16-1,Garijo14}]\label{conj}
If $G$ is an isolate-free and twin-free graph of order $n$, then  $\gamma_{LD}(G)\le \frac{n}{2}$.
\end{conj}

Garijo et al.~\cite{Garijo14} showed that if $G $ is a twin-free graph of order $n$, then $\gamma_{LD}(G)\le \lfloor\frac{2n}{3}\rfloor+1$. Later, Foucaud et al.~\cite{Foucaud16-2} subsequently improved this upper bound to $\lfloor\frac{2n}{3}\rfloor$. Recently, Bousquet et al.~\cite{Bousquet25} further reduced the upper bound to $\lceil\frac{5n}{8}\rceil$, which is currently the best known bound to Conjecture~\ref{conj}.  Conjecture~\ref{conj} has not yet been proven, but has been shown to be true for some important graph classes. 

\begin{theorem}\label{supportconj}
Conjecture~\ref{conj} is true for isolate-free and twin-free graph $G$ of orders $n$ if at least one of the following is satisfied.
\begin{enumerate}
	\item \rm{\cite{Garijo14}} $G$ has no $4$-cycle.
	\item \rm{\cite{Garijo14}} $G$ has independence number at least $\frac{n}{2}$.
	\item \rm{\cite{Garijo14}} $G$ has clique number at least $\lceil\frac{n}{2}\rceil+1$.
	\item \rm{\cite{Balbuena15}} $G$ has girth at least $5$ and minimum degree at least $2$.
	\item \rm{\cite{Foucaud16-2}} $G$  is a split graph or a co-bipartite graph.
	\item \rm{\cite{Foucaud17}} $G$ is a line graph.
	\item \rm{\cite{Claverol21}} $G$ is a maximal outerplanar graph.
	\item \rm{\cite{Chakraborty24}} $G$ is a block graph.
	\item \rm{\cite{Chakraborty24-1}} $G$ is a subcubic graph.
\end{enumerate}
\end{theorem}

Given a graph $G$, if there exist two LD-sets $D_1$ and $D_2$ such that $D_1\cup D_2=V(G)$, $D_1\cap D_2=\emptyset$, then $[D_1,D_2]$ is called an \emph{LD-partition} of $G$. Motivated by Conjecture~\ref{conj}, several authors have explored the following, a slightly stronger question.

\begin{quest}[\cite{Chakraborty25,Foucaud16-1, Foucaud16-2, Garijo14}]\label{motivation}
	For an isolate-free and twin-free graph $G$, does $G$ admit  an LD-partition$?$
\end{quest}
  
  Recently, Chakraborty et al.~\cite{Chakraborty25}
  showed that if $G$ is an isolate-free (and not necessarily twin-free) graph, then $V(G)$ can be partitioned into a dominating set and an LD-set. It is already known that Question~\ref{motivation} has a positive answer for bipartite graphs~\cite{Garijo14} and block graphs~\cite{Chakraborty24}, which naturally motivates investigating its validity in superclasses of these graph classes. In this context, we provide a positive answer to Question~\ref{motivation} for the class of distance-hereditary graphs, a well-known superclass of block graphs. Additionally, while Conjecture~\ref{conj} is known to hold for maximal outerplanar graphs, split graphs, and co-bipartite graphs, we go further by establishing an affirmative answer to Question~\ref{motivation} for each of these graph classes as well. 

   This paper is organized as follows. In Section~\ref{sec:D_H}, we show that if $G$ is an isolate-free and twin-free distance-hereditary graph, then $G$ admits an LD-partition. In Section~\ref{sec:mop}, we prove that every maximal outerplanar graph of order at least 4 admits an LD-partition. In Section~\ref{sec:spco}, we prove that every isolate-free and twin-free split graph and co-bipartite graph also admits an LD-partition. Finally, in Section~\ref{sec:con}, we discuss potential directions for future research.

\section{Distance-hereditary graphs}\label{sec:D_H}

Our objective in this section is to show that every twin-free and isolate-free distance-hereditary graph admits an LD-partition. In this section, we assume that $G$ is a connected distance-hereditary graph. 

\medskip
Chang et al.~\cite{Chang97} characterized distance-hereditary graphs via edge connections between two special sets of vertices, called twin sets. The comprehensive procedure is given in the next paragraph. At its base level, a graph $G$ with a single vertex $v$ is recognized as a distance-hereditary graph, endowed with the twin set $TS(G) = \{v\}$.

\medskip
A distance-hereditary graph $G$ can be constructed from two existing vertex-disjoint distance-hereditary graphs, $G_l$ and $G_r$, each possessing twin sets $TS(G_l)$ and $TS(G_r)$, respectively, by using any of the subsequent three operations.

\begin{itemize}
    
    \item If the \emph{true twin} operation $\otimes$ is applied to construct the graph $G$ from $G_l$ and $G_r$, then
    \begin{itemize}
        \item The vertex set of $G$ is  $V(G) = V(G_l) \cup V(G_r)$.
        \item The edge set of $G$ is  $E(G) = E(G_l) \cup E(G_r) \cup \{v_1v_2 \vert v_1 \in TS(G_l), v_2 \in TS(G_r)\}$.
        \item The twin set of $G$ is  $TS(G) = TS(G_l) \cup TS(G_r)$.
    \end{itemize}
    \item If the \emph{false twin} operation $\odot$ is employed to construct the graph $G$ from $G_l$ and $G_r$, then
    \begin{itemize}
        \item The vertex set of $G$ is  $V(G) = V(G_l) \cup V(G_r)$.
        \item The edge set of $G$ is  $E(G) = E(G_l) \cup E(G_r)$.
        \item The twin set of $G$ is  $TS(G) = TS(G_l) \cup TS(G_r)$.
    \end{itemize}
    \item If the \emph{attachment} operation $\oplus$ is employed to construct the graph $G$ from $G_l$ and $G_r$, then
    \begin{itemize}
        \item The vertex set of $G$ is $V(G) = V(G_l) \cup V(G_r)$.
        \item The edge set of $G$ is $E(G) = E(G_l) \cup E(G_r) \cup \{v_1v_2 ~\vert~v_1 \in TS(G_l), v_2 \in TS(G_r)\}$.
        \item The twin set of $G$ is $TS(G) = TS(G_l)$.
    \end{itemize}
\end{itemize}

\begin{figure}[ht]
\centering
\begin{subfigure}{.45\textwidth}
  \centering
  % include first image
  \includegraphics[width=.7\linewidth]{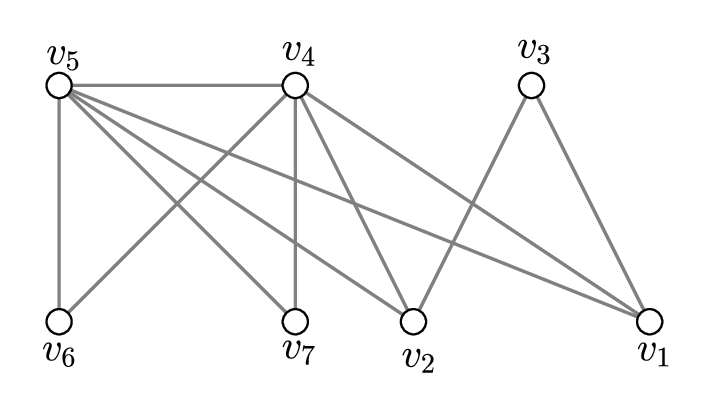}  
  \caption{A distance-hereditary graph $G$}
\end{subfigure}
\begin{subfigure}{.4\textwidth}
  \centering
  % include second image
  \includegraphics[width=0.9\linewidth]{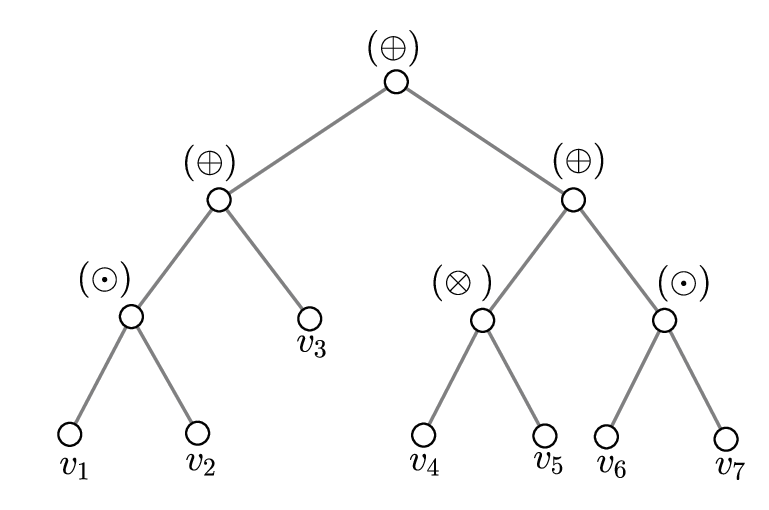}  
  \caption{The decomposition tree $T_G$ of $G$}
\end{subfigure}
\caption{An example of a distance-hereditary graph with its decomposition tree}
\label{fig:3}
\end{figure}

\medskip
By employing the three operations detailed above, one can systematically construct any distance-hereditary graph. This process leads to the creation of a binary tree representation for a given distance-hereditary graph $G$, commonly referred to as a \emph{decomposition tree}. The definition of this tree is structured as follows: it articulates the sequence of operations through a full binary tree $T$, where the leaves of $T$ correspond to the vertices of $G$. Furthermore, each internal vertex in $T$ is assigned one of the labels $\otimes, \odot,$ or $\oplus$, signifying the true twin operation, false twin operation, and attachment operation, respectively.

\medskip
In this representation, each leaf of $T$ corresponds to a distance-hereditary graph with a single vertex. A rooted subtree $T'$ of $T$ corresponds to the induced subgraph of $G$ on the vertices represented by the leaves of $T'$. Note that this induced subgraph is itself a distance-hereditary graph. For an internal vertex $v$ of $T$, the label of $v$ corresponds to the operation between the subgraphs represented by the subtrees rooted at the left and right children of $v$. Note that the order of the children only matters for the $\oplus$ operation. An example is illustrated in Figure \ref{fig:3}.

\medskip
 Next, we prove the main theorem of this section.

\begin{theorem}\label{th:main_theorem_D_H}

If $G$ is an isolate-free and twin-free distance-hereditary graph, then $G$ admits an LD-partition.
\end{theorem}

\begin{proof}

We prove this using induction on $\vert V(G)\vert+ \vert E(G)\vert=n+m$. For the base cases, it can be easily checked for all isolate-free and twin-free distance-hereditary graphs of order $4$ (as the only example of twin-free and isolate-free distance-hereditary graph of order $4$ is $P_4$). So, let the statement be true for all isolate-free and twin-free distance-hereditary graphs of order $<n+m$.

    \medskip
    Let $G$ be an isolate-free and twin-free distance-hereditary graph of order $n$. Suppose $G$ is disconnected. Let $G_1$,$G_2$,\ldots,$G_k$ be the components of $G$. Then by the induction hypothesis, each $G_i$ admits an LD-partion $[D^{i}_1,D^{i}_2]$ for all $i\in\{1,2,\ldots,k\}$. Let $D_1=D^{1}_1\cup D^{2}_1\cup \ldots \cup D^{k}_1$ and $D_2=D^{1}_2\cup D^{2}_2\cup \ldots \cup D^{k}_2$. Then $[D_1,D_2]$ is an LD-partition of $G$. So we assume that $G$ is connecetd. Let $T_G$ be the decomposition tree of $G$. Consider the BFS levels of $T_G$ and let $t$ be an internal node of $T_G$ that is situated at the second last BFS level. Note that both children (say $a$ and $b$) of $t$ are vertices of $G$.

    Suppose $t$ has label $\odot$ (or $\otimes$). This implies that $a$ and $b$ are false twins (or true twins), which contradicts the fact that $G$ is twin-free. Hence, $t$ has label $\oplus$. Let $t'$ be the parent of $t$. Now, a vertex $t'$ can be labeled as $\odot$, $\otimes$, or $\oplus$. Moreover, if $t'$ has the label $\oplus$ and its child is a leaf node, then two cases arise: either $t$ is the left child of $t'$, or $t$ is the right child of $t'$. Moreover, if the child of $t'$ other than $t$, say $t''$, is an internal node, then $t''$ has the label $\oplus$ (a similar argument can be given as the one used to show that $t'$ has the label $\oplus$). Based on these observations, we consider the following cases.
    
    \begin{enumerate}
    	\item [(1)] $t'$ has label $\oplus$, $t$ is the right child of $t'$, and the left child of $t'$ is a leaf node $c$;
    	\item [(2)] $t'$ has label $\oplus$, $t$ is the left child of $t'$, and the right child of $t'$ is a leaf node $c$;
    	\item [(3)] $t'$ has label $\odot$ and the other child of $t'$ is a leaf node $c$; 
    	\item [(4)] $t'$ has label $\otimes$ and the other child of $t'$ is a leaf node $c$; 
    	\item [(5)] $t'$ has label $\oplus$ and the other child of $t'$ is an internal node $t''$ which also has label $\oplus$; 
    	\item [(6)] $t'$ has label $\odot$ and the other child of $t'$ is an internal node $t''$ which has label $\oplus$; 
    	\item [(7)] $t'$ has label $\otimes$ and the other child of $t'$ is an internal node $t''$ which has label $\oplus$. 
    \end{enumerate}
    We now analyze these cases separately. 

\medskip
\noindent \textbf{Case 1}: $t'$ has label $\oplus$, $t$ is the right child of $t'$, and the left child of $t'$ is a node $c$. For a clear understanding, see Figure \ref{fig:DH_1}.

\begin{figure}[ht]
    \centering
\begin{subfigure}{.45\textwidth}
  \centering
  % include first image
  \includegraphics[width=.4\linewidth]{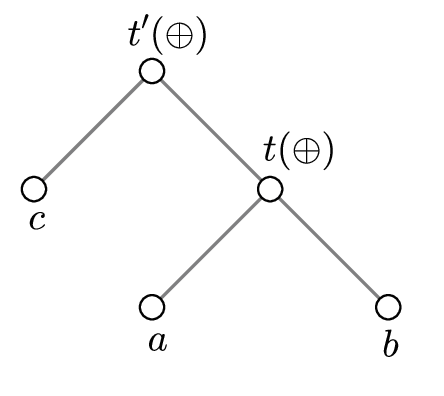}  
  \caption{Subtree of $T_G$ rooted at $t'$}
  
\end{subfigure}
\begin{subfigure}{.35\textwidth}
  \centering
  % include second image
  \includegraphics[width=1\linewidth]{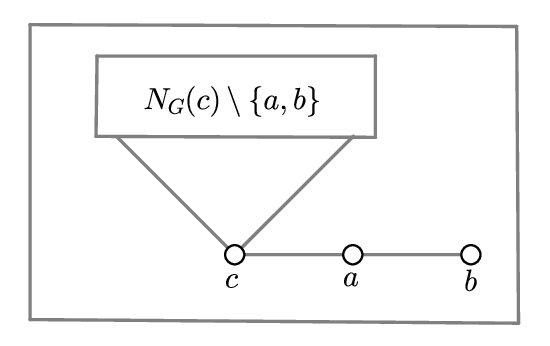}  
  \caption{$G$}
\end{subfigure}
    \caption{The subtree of $T_G$ rooted at $t'$ in Case 1}
    \label{fig:DH_1}
\end{figure}

\medskip
Let $G_1=G\setminus \{a,b\}$. If $G_1$ is twin-free and isolate-free, then by the induction hypothesis, $G_1$ admits an LD-partition $[D'_1,D'_2]$. Without loss of generality, let $c\in D'_1$ and $c\notin D'_2$. Define $D_1=D'_1\cup\{b\}$ and $D_2=D'_2\cup \{a\}$. Observe that $[D_1,D_2]$ is an LD-partition of $G$. If $G_1$ is not twin-free, then there exists a vertex in $V(G_1)$, say $x$, such that $c$ and $x$ are twins in $G_1$. In the following, we consider two cases and in each case, we prove that $G$ admits an LD-partition.

\medskip
\noindent \textbf{Case 1.1}: $c$ and $x$ are false twins in $G_1$.

\begin{figure}[ht]
\centering
\begin{subfigure}{.4\textwidth}
  \centering
  % include first image
  \includegraphics[width=.7\linewidth]{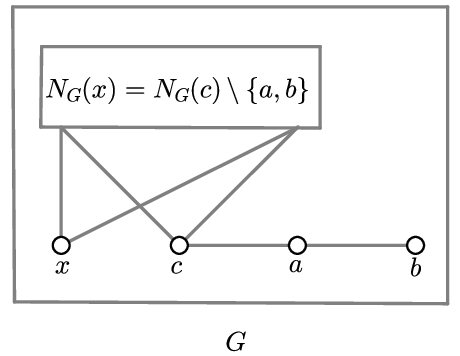}  
  \caption{Case 1.1}
 
\end{subfigure}
\begin{subfigure}{.4\textwidth}
  \centering
  % include second image
  \includegraphics[width=0.7\linewidth]{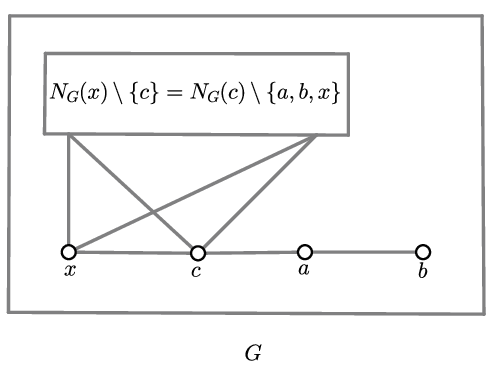}  
  \caption{Case 1.2}
  
\end{subfigure}
\caption{Case 1}
\label{fig:case1_1_2}
\end{figure}

\medskip
Let $G'=G\setminus\{a,b,c\}$.  Suppose there exist twins in $G'$, say $y$ and $z$. Clearly, exactly one of $y$ and $z$ is adjacent to $c$, which implies that exactly one of $y$ and $z$ is adjacent to $x$. Since $x\in V(G')$, $N_{G'}(y)\neq N_{G'}(z)$, contradicting that $y$ and $z$ are twins. Hence, $G'$ is twin-free. Then by the induction hypothesis, let $[D'_1,D'_2]$ be an LD-partition of $G'$ such that $x\in D'_1$ and $x\notin D'_2$. We define $D_1=D'_1\cup \{c,b\}$ and $D_2=D'_2\cup \{a\}$. It is easy to check that $D_1$ is an LD-set of $G$.

\medskip
Clearly, $D_2$ is a dominating set of $G$. For the sake of contradiction, assume that $D_2$ is not an LD-set of $G$. This implies that there exist two vertices in $V(G)\setminus D_2$ which have the same neighborhood in $D_2$. The only candidates for these two vertices are $c$ and $b$ (as the rest of the vertices in $V(G)\setminus D_2$ have different neighborhoods in $D_2$ since $D'_2$ is an LD-set of $G'$). But note that $x$ is dominated by some vertex $z$ in $D'_2$, so $N_G(z)$ contains $c$ but not $b$, which contradicts the fact that $c$ and $b$ have the same neighborhood in $D_2$. This implies that $D_2$ is an LD-set. So $[D_1,D_2]$ is an LD-partition of $G$.

\medskip
\noindent \textbf{Case 1.2}: $c$ and $x$ are true twins of $G_1$.

\medskip
 Let $G'=G\setminus \{b\}$. Note that $G'$ is a twin-free distance-hereditary graph. By the induction hypothesis, let $[D'_1,D'_2]$ be an LD-partition of $G'$. If $\{c,a\}\subseteq D'_1$, then $c$ and $a$ do not belong to $D'_2$, which contradicts the fact that $D'_2$ is a dominating set. So $\{c,a\}\not\subseteq D'_1$. Similarly, $\{c,a\}\not\subseteq D'_2$. Hence, without loss of generality, let $c\in D'_1$ and $a\in D'_2$. We define $D_1=D'_1\cup \{b\}$ and $D_2=D'_2$. Clearly, $D_1$ is an LD-set of $G$ and $D_2$ is a dominating set of $G$. For the sake of contradiction, assume that $D_2$ is not an LD-set of $G$, this implies that $N_G(c)\cap D_2=N_G(b)\cap D_2=\{a\}$. This means $N_{G'}[c] \cap D_2=N_{G'}[c] \cap D'_2=\{a\}$, implying $N_{G'}[x]\cap D'_2=\emptyset$ which contradicts the fact that $D'_2$ is a dominating set of $G'$. Hence, $[D_1,D_2]$ is an LD-partition of $G$.

\medskip
\noindent \textbf{Case 2}: $t'$ has label $\oplus$, $t$ is the left child of $t'$, and the right child of $t'$ is a node $c$. For clear understanding, see Figure \ref{fig:case2}.

\begin{figure}[ht]
    \centering
    \begin{subfigure}{.45\textwidth}
  \centering
  % include first image
  \includegraphics[width=.4\linewidth]{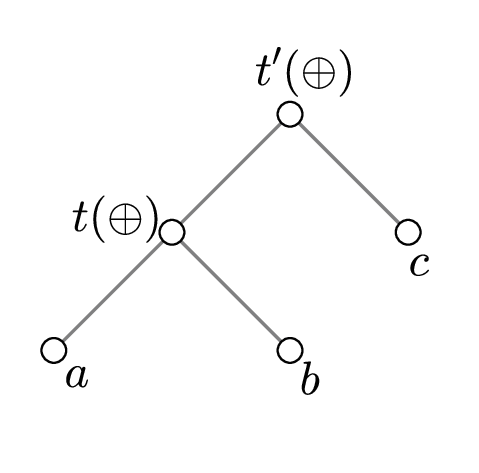}  
  \caption{Subtree of $T_G$ rooted at $t'$}
\end{subfigure}
\begin{subfigure}{.3\textwidth}
  \centering
  % include second image
  \includegraphics[width=0.55\linewidth]{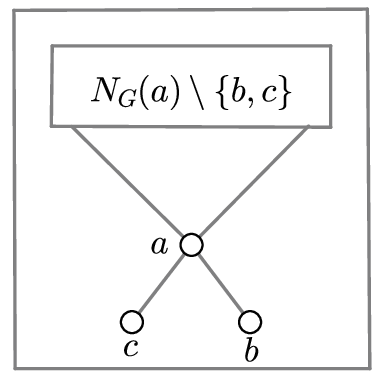}  
  \caption{$G$}
\end{subfigure}
    \caption{The subtree of $T_G$ rooted at $t'$ in Case 2}
    \label{fig:case2}
\end{figure}

\medskip
In this case, $b$ and $c$ are leaves in $G$ that are adjacent to the vertex $a$. Hence $b$ and $c$ are twins, which contradicts the fact that $G$ is twin-free. So this is not a valid case.

\medskip
\noindent \textbf{Case 3}: $t'$ has label $\odot$ and the other child of $t'$ is a leaf node $c$. For clear understanding, see Figure \ref{fig:case3}.

\begin{figure}[ht]
    \centering
    \begin{subfigure}{.45\textwidth}
  \centering
  % include first image
  \includegraphics[width=.45\linewidth]{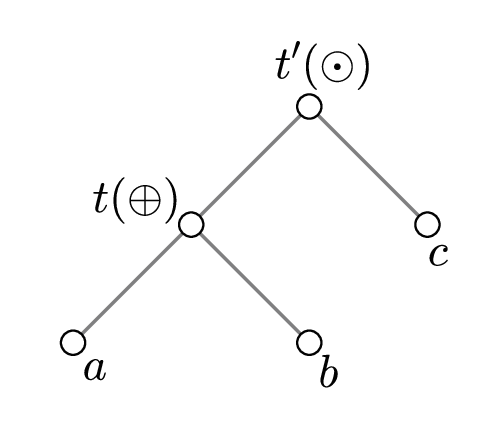}  
  \caption{Subtree of $T_G$ rooted at $t'$}
  
\end{subfigure}
\begin{subfigure}{.3\textwidth}
  \centering
  % include second image
  \includegraphics[width=0.65\linewidth]{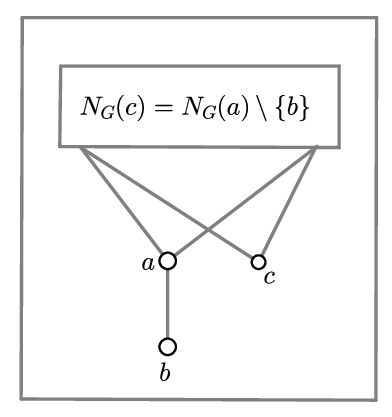}  
  \caption{$G$}
\end{subfigure}
    \caption{The subtree of $T_G$ rooted at $t'$ in Case 3}
    \label{fig:case3}
\end{figure}

\medskip
Let $G'=G\setminus \{a,b\}$. By similar arguments as we did at the beginning of Case 1.1, it can be shown that $G'$ is also twin-free and isolate-free. By the induction hypothesis, let $[D'_1,D'_2]$ be an LD-partition of $G'$. Without loss of generality, let $c\in D'_1$ and $c\notin D'_2$.  We define $D_1=D'_1\cup \{a\}$ and $D_2=D'_2\cup \{b\}$. By using analogous arguments like in Case 1.1, it can be shown that $[D_1,D_2]$ is an LD-partition of $G$.

\medskip
\noindent \textbf{Case 4}: $t'$ has label $\otimes$ and the other child of $t'$ is a leaf node $c$. For clear understanding, see Figure \ref{fig:case4}.

\medskip
Let $G'$ be obtained from $G$ by deleting all the edges between $N_G(c)\setminus \{a\}$ and $a$. By similar arguments as we did at the beginning of Case 1.1, it can be proved that $G'$ is twin-free. Let $G''=G'\setminus \{a,b\}$. Note that $G'$ and
$G''$ are isolate-free. In the following, we consider two cases and in each case, we prove that $G$ admits an LD-partition.

\begin{figure}[ht]
    \centering
    \begin{subfigure}{.45\textwidth}
  \centering
  % include first image
  \includegraphics[width=.4\linewidth]{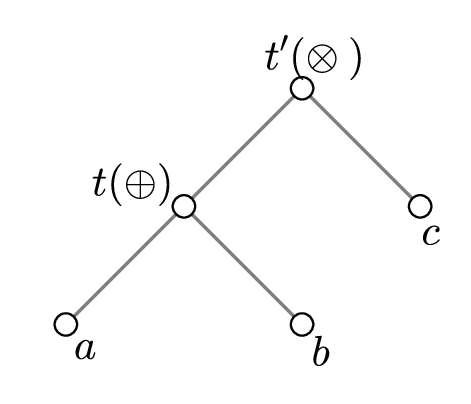}  
  \caption{Subtree of $T_G$ rooted at $t'$}
  
\end{subfigure}
\begin{subfigure}{.3\textwidth}
  \centering
  % include second image
  \includegraphics[width=0.7\linewidth]{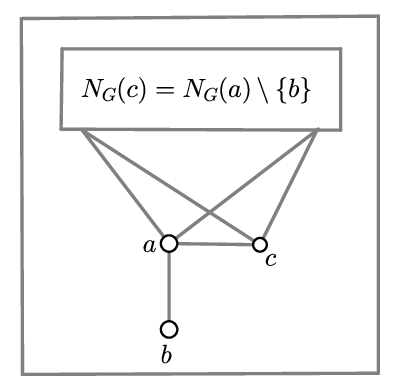}  
  \caption{$G$}
\end{subfigure}
    \caption{The subtree of $T_G$ rooted at $t'$ in Case 4}
    \label{fig:case4}
\end{figure}

\medskip
\noindent \textbf{Case 4.1}: $G''$ is twin-free.

\medskip
By the induction hypothesis, $G''$ admits an LD-partition $[D'_1,D'_2]$ such that $c\in D'_1$.  We define $D_1=D'_1\cup \{a\}$ and $D_2=D'_2\cup \{b\}$. By using analogous arguments as in Case 1.1, it can be shown that $[D_1,D_2]$ is an LD-partition of $G$.

\medskip
\noindent \textbf{Case 4.2}: $G''$ is not twin-free. 

\begin{figure}[ht]
    \centering
    \includegraphics[width=0.3\linewidth]{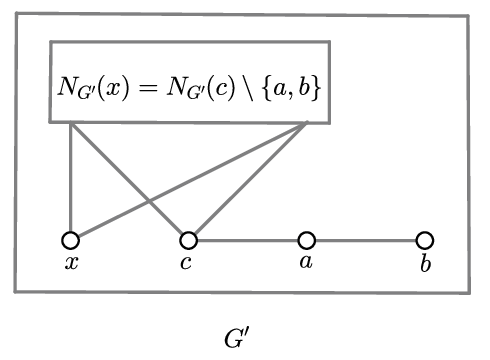}
    \caption{Case 4.2 ($c$ and $x$ are false twins in $G''$)}
    \label{fig:case4_2}
\end{figure}

Since $G''$ is not twin-free, there exists a vertex in $V(G'')$, say $x$, such that $c$ and $x$ are twins in $G''$. First we prove that $c$ and $x$ are false twins. For the sake of contradiction, let $c$ and $x$ be true twins in $G''$, implying $x\in N_{G}(c)\subseteq N_{G}(a)$. Hence $N_{G''}(x)=N_{G''}(c)$ implies that  $N_{G}[x]=N_{G}[c]$, which contradicts the fact that $G$ is twin-free. Hence, $c$ and $x$ must be false twins in $G''$ (refer to Figure \ref{fig:case4_2}).

Recall that $G'$ is twin-free. Note that $G'$ has same structure as in Case 1.1. Hence, by the analysis of Case 1.1 and induction hypothesis, it can be concluded that, $G'$ admits an LD-partition $[D_1,D_2]$ such that $D_1$ contains $x,c,b$ and $D_2$ contains $a$. It is easy to see that $D_1$ is also an LD-set of $G$.

Next, we prove that $D_2$ is an LD-set of $G$. Note that $D_2$ is a dominating set of $G$, as it is a dominating set of $G'$. For the sake of contradiction, assume that $D_2$ is not an LD-set of $G$. This implies that adding back the deleted edges has created some problem. Hence, there exists $y\in V(G)$ such that $N_{G}(y)\cap D_2=N_G(c)\cap D_2$. This implies that $N_{G'}(y)\cap D_2=(N_{G'}(c)\cap D_2)\setminus \{a\}$, implying $N_{G'}(y)\cap D_2=N_{G'}(x)\cap D_2$ which contradicts the fact that $D_2$ is an LD-set in $G'$. Hence $D_2$ is an LD-set in $G$. So $[D_1,D_2]$ is an LD-partition of $G$.

\medskip
\noindent \textbf{Case 5}: $t'$ has label $\oplus$ and the other child of $t'$ is an internal node $t''$ which also has label $\oplus$. For a clear understanding, see Figure \ref{fig:case5}.

\begin{figure}[ht]
    \centering
     \begin{subfigure}{.5\textwidth}
  \centering
  % include first image
  \includegraphics[width=.4\linewidth]{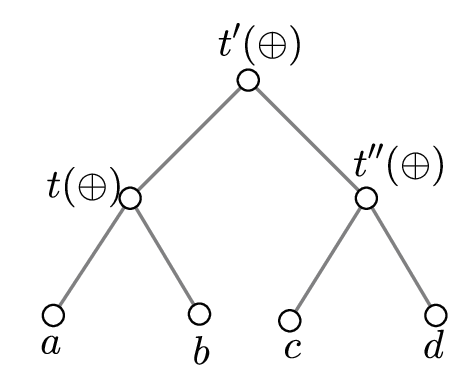}  
  \caption{Subtree of $T_G$ rooted at $t'$}
  
\end{subfigure}
\begin{subfigure}{.4\textwidth}
  \centering
  % include second image
  \includegraphics[width=0.45\linewidth]{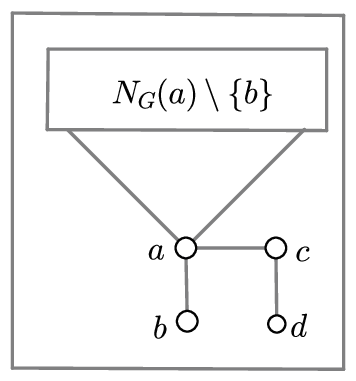}  
  \caption{$G$}
\end{subfigure}
    \caption{Case 5}
    \label{fig:case5}
\end{figure}

\medskip
Let $G'=G\setminus \{c,d\}$. Observe that $G'$ is a twin-free and isolate-free distance-hereditary graph. By the induction hypothesis, $G'$ admits an LD-partition $[D'_1,D'_2]$ such that $D'_1$ contains $a$ and $D'_2$ contains $b$. We define $D_1=D'_1\cup \{c\}$ and $D_2=D'_2\cup \{d\}$. It is easy to observe that $[D_1,D_2]$ is an LD-partition of $G$.

\medskip
\noindent \textbf{Case 6}: $t'$ has label $\odot$ and the other child of $t'$ is an internal node $t''$ which has label $\oplus$. Refer to Figure \ref{fig:case6} for clear understanding.

\begin{figure}[ht]
    \centering
    \begin{subfigure}{.5\textwidth}
  \centering
  % include first image
  \includegraphics[width=.4\linewidth]{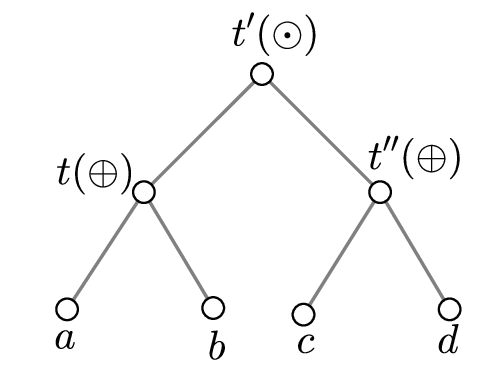}  
  \caption{Subtree of $T_G$ rooted at $t'$}
  \label{fig:sub1}
\end{subfigure}
\begin{subfigure}{.4\textwidth}
  \centering
  % include second image
  \includegraphics[width=0.55\linewidth]{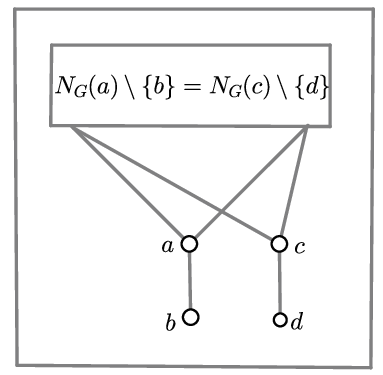}  
  \caption{$G$}
  \label{fig:sub2}
\end{subfigure}
    \caption{Case 6}
    \label{fig:case6}
\end{figure}

Let $G'=G\setminus \{c,d\}$. Observe that $G'$ is a twin-free and isolate-free distance-hereditary graph. By the induction hypothesis, $G'$ admits an LD-partition $[D'_1,D'_2]$ such that $D'_1$ contains $a$ and $D'_2$ contains $b$. We define $D_1=D'_1\cup \{c\}$ and $D_2=D'_2\cup \{d\}$. It is easy to observe that $[D_1,D_2]$ is an LD-partition of $G$.

\medskip
\noindent \textbf{Case 7}: $t'$ has label $\otimes$ and the other child of $t'$ is an internal node $t''$ which has label $\oplus$. Refer to Figure \ref{fig:case7} for clear understanding.

\begin{figure}[ht]
    \centering
    \begin{subfigure}{.5\textwidth}
  \centering
  % include first image
  \includegraphics[width=.37\linewidth]{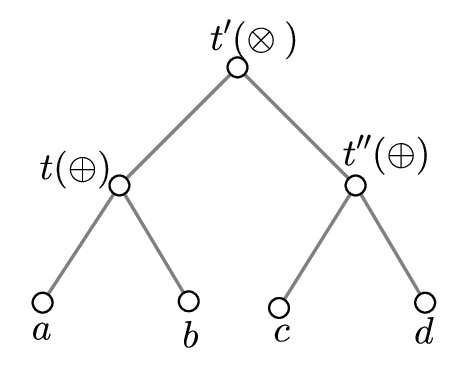}  
  \caption{Subtree of $T_G$ rooted at $t'$}
\end{subfigure}
\begin{subfigure}{.3\textwidth}
  \centering
  % include second image
  \includegraphics[width=0.8\linewidth]{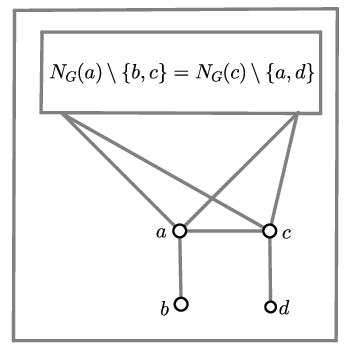}  
  \caption{$G$}
\end{subfigure}
    \caption{Case 7}
    \label{fig:case7}
\end{figure}

Let $G'=G\setminus \{c,d\}$. Observe that $G'$ is a twin-free and isolate-free distance-hereditary graph. By the induction hypothesis, $G'$ admits an LD-partition $[D'_1,D'_2]$ such that $D'_1$ contains $a$ and $D'_2$ contains $b$. We define $D_1=D'_1\cup \{c\}$ and $D_2=D'_2\cup \{d\}$. It is easy to observe that $[D_1,D_2]$ is an LD-partition of $G$.

This completes the proof of Theorem \ref{th:main_theorem_D_H}.\end{proof}

\section{Maximal Outerplanar graphs}\label{sec:mop}

 Note that the only mops with twin vertices are those of order 3 and 4: the cycle of order 3, and the graph obtained by deleting an edge from the complete graph of order 4. In question~\ref{motivation}, we are only interested in the twin-free graphs. However, in this section, we prove that every mop of order at least $4$ has an LD-partition.  To prove the above result, we need the following observation.

\begin{obs}
	\label{obs1}
	If $G$ is a mop of order $4$ or $5$, then there exists a vertex adjacent to all other vertices of $G$.
\end{obs}

Let $G$ be a mop of order $n \ge 4$ vertices. Hence, there exists a plane embedding of $G$ such that all vertices of $G$ are on the outer face, and all inner faces are triangles. We construct a new graph $T$ associated with a given mop $G$ as follows.\\ [-20pt]
\begin{itemize}
\item   Each vertex of $T$ represents a triangle in $G$.\\ [-20pt]
\item  Two vertices in $T$ are adjacent by an edge if their corresponding triangles in $G$ share an edge.\\ [-20pt]
\end{itemize}
 Note that $T$ is connected. If $T$ has a cycle, then there exists a vertex in $G$ that is enclosed by triangles, which is not possible since $G$ is outerplanar. So $T$ must be a tree. The maximum degree of any vertex in $T$ is at most 3. We will analyze the tree $T$ and understand its corresponding structure in mop $G$.
\begin{theorem}\label{thmmop}
	If $G$ is a mop of order at least $4$, then $G$ admits an LD-partition.
\end{theorem}
    
\begin{proof}
	Let $G$ be a mop of order $n\ge 4$. We will use induction on $n$. Let $v_1v_2v_3\ldots v_nv_1$ be the vertices on the outer face of $G$, listed in order. If $n=4$, then by Observation~\ref{obs1}, without loss of generality, assume that $v_1$ is adjacent to $v_2,v_3$, and $v_4$ in $G$. Since $G$ is a mop, $v_2v_4\notin E(G)$. Then $\{v_1,v_2\}$ and $\{v_3,v_4\}$ are LD-sets of $G$. If $n=5$, then by Observation~\ref{obs1}, without loss of generality, assume that $v_1$ is adjacent to $v_2,v_3,v_4$, and $v_5$ in $G$. Since $G$ is a mop, $v_2v_4, v_2v_5, v_3v_5\notin E(G)$. Then $\{v_1,v_2,v_5\}$ and $\{v_3,v_4\}$ are LD-sets of $G$. So assume that $n\ge 6$.

    \medskip
	Let $T$ be the tree associated with the mop $G$, where $T$ is rooted at a leaf $w$. Since $T$ has at least two leaves, there exists a leaf other than $w$, say $y$. Since $n\ge 6$, $|V(T)|\ge 4$. Let $x$ be the parent of $y$ in $T$. We define $T_x$ as the subtree of $T$ that is rooted at the vertex $x$. Since the maximum degree of $T$ is at most 3, $x$ has at most two children in $T$. In the following, we consider two cases based on the number of children of $x$.
		\medskip 
		
	\noindent\textbf{Case 1:} $x$ has two children in $T$.
	
	\medskip 
	
	 Let $z$ be a child of $x$ other than $y$. Let $R_x$ be the triangle in $G$ corresponding to the vertex $x$. Let $V(R_x) = \{v_1,v_2,v_3\}$. Let $R_y$ and $R_z$ be the triangles in $G$ corresponding to the vertices $y$ and $z$, respectively. Further, let $V(R_y) = \{v_1,v_3,v_4\}$ and $V(R_z) = \{v_2,v_3,v_5\}$. Thus, $G$ contains the subgraph illustrated in Figure~\ref{fig:paras_1}(b), where the shaded triangle corresponds to the vertex $x$ in~$T_x$. Since $y$ and $z$ are leaves in $T$, we have that $\deg_G(v_4) = \deg_G(v_5) = 2$ and $\deg_G(v_3) = 4$. Recall that $n \ge 6$.
		
		\begin{figure}[ht]
        \centering
			\begin{subfigure}{.4\textwidth}
				\centering
				% include first image
				\includegraphics[width=.25\linewidth]{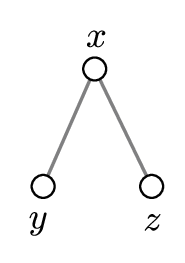}  
				\caption{Subtree $T_x$ of $T$}
				\label{fig:paras_sub1}
			\end{subfigure}
			\begin{subfigure}{.4\textwidth}
				\centering
				% include second image
				\includegraphics[width=.5\linewidth]{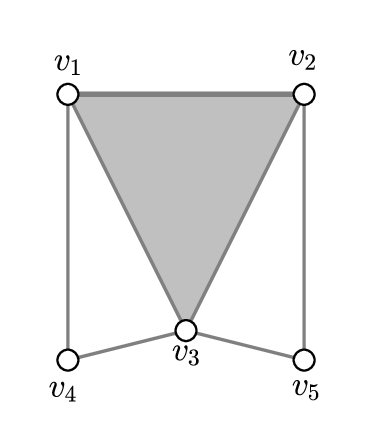}  
				\caption{Subgraph of $G$ corresponding to subtree $T_x$}
				\label{fig:paras_sub2}
			\end{subfigure}
			\caption{Subtree $T_x$ and possible subgraph of $G$.}
			\label{fig:paras_1}
		\end{figure}
		\medskip
		Let $H$ be the graph of order $n'$ obtained from $G$ by deleting the vertices $v_4$ and $v_5$. Since $n \ge 6$, we have $n' = n - 2 \ge 4$. We note that $H$ is also a mop. Then by the induction hypothesis, let $[D'_1,D'_2]$ be an LD-partition of $H$. Since $\deg_{H}(v_3)=2$, each $D'_i$ contains at least one vertex from the set $\{v_1,v_2,v_3\}$ for $i\in\{1,2\}$. In the following, we consider three cases and in each case, we prove that $G$ admits an LD-partition.
	
    \medskip	
		\noindent\textbf{Case 1.1:} $v_1,v_2\in D'_1$ and $v_3\in D'_2$.

        \medskip 
		Let $D_1=D'_1\cup\{v_3\}$ and $D_2=(D'_2\setminus\{v_3\})\cup\{v_4,v_5\}$. Now we show that each $D_i$ is an LD-set of $G$ for $i\in\{1,2\}$. Since $N_G(v_4)\cap D_1=\{v_1,v_3\}$ and $N_G(v_5)\cap D_1=\{v_2,v_3\}$, $D_1$ is an LD-set of $G$. Note that $N_G(v_1)\cap\{v_4,v_5\}=\{v_4\}$, $N_G(v_2)\cap\{v_4,v_5\}=\{v_5\}$, and $N_G(v_3)\cap\{v_4,v_5\}=\{v_4,v_5\}$. Moreover, $\{v_4,v_5\}\subset D_2$. Therefore, $D_2$ is an LD-set of $G$.
	
    \medskip			
	   \noindent\textbf{Case 1.2:} $v_1,v_3\in D'_1$ and $v_2\in D'_2$.
       
	\medskip    
	    Let $D_1=D'_1\cup\{v_5\}$ and $D_2=D'_2\cup\{v_4\}$. Now we show that each $D_i$ is an LD-set of $G$ for $i\in\{1,2\}$. Note that $N_G(v_2)\cap\{v_1,v_3,v_5\}=\{v_1,v_3,v_5\}$ and  $N_G(v_4)\cap\{v_1,v_3,v_5\}=\{v_1,v_3\}$. Moreover, $\{v_1,v_3,v_5\}\subset D_1$. Therefore, $D_1$ is an LD-set of $G$. Since $D'_2$ is an LD-set of $H$, we have $N_H(v_1)\cap D'_2\neq N_H(v_3)\cap D'_2$. Hence $N_G(v_1)\cap D_2\neq N_G(v_3)\cap D_2$. Note that $\{v_2,v_4\}\subseteq N_G(v_1)\cap D_2$, $N_G(v_3)\cap D_2= \{v_2,v_4\}$, and  $N_G(v_5)\cap D_2=\{v_2\}$. Therefore, $D_2$ is an LD-set of $G$.
	 
     \medskip  
     
		\noindent\textbf{Case 1.3:} $v_2,v_3\in D'_1$ and $v_1\in D'_2$.
        
	\medskip 	
		Let $D_1=D'_1\cup\{v_4\}$ and $D_2=D'_2\cup\{v_5\}$.  We can make similar arguments as we did in Case 1.2.
		
		Hence $G$ admits an LD-partition $[D_1,D_2]$. 
    	\medskip 

	\noindent\textbf{Case 2:} $x$ has only one child in $T$.
	
   \medskip 
   	
    Since $n\ge 6$, $|V(T)|\ge 4$. So the parent of $x$ exists in $T$. Let $x'$ be the parent of $x$ in $T$. Recall that $x$ is the parent of $y$ in  $T$. Let $R_{x'}$ be the triangle in $G$ corresponding to the vertex $x'$. Let $V(R_{x'}) = \{v_1,v_2,v_3\}$. Let $R_x$ and $R_y$ be the triangles in $G$ corresponding to the vertices $x$ and $y$, respectively.  Further, due to symmetry, without loss of generality, assume that $V(R_x) = \{v_1,v_2,v_4\}$ and $V(R_y) = \{v_2,v_4,v_5\}$. Thus, $G$ contains the subgraph illustrated in Figure~\ref{fig:paras_2}(b), where the shaded triangle corresponds to the vertex $x'$ in~$T_{x'}$. Since $y$ is a leaf and $x$ has exactly one child in $T$, we have $\deg_G(v_5) = 2$ and $\deg_G(v_4) =3$.

\begin{figure}[ht]
	\centering
	\begin{subfigure}{.4\textwidth}
		\centering
		% include first image
		\includegraphics[width=.15\linewidth]{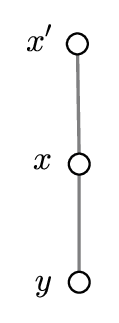}  
		\caption{Subtree $T_{x'}$ of $T$}
		\label{fig:paras_sub3}
	\end{subfigure}
	\begin{subfigure}{.4\textwidth}
		\centering
		% include second image
		\includegraphics[width=.6\linewidth]{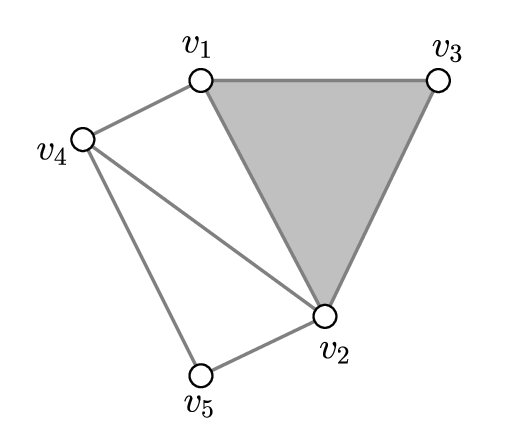}  
		\caption{Subgraph of $G$ corresponding to subtree $T_{x'}$}
		\label{fig:paras_sub4}
	\end{subfigure}
	\caption{Subtree $T_{x'}$ and possible subgraph of $G$.}
	\label{fig:paras_2}
\end{figure}

\medskip 
Let $H$ be the graph of order $n'$ obtained from $G$ by deleting the vertices $v_4$ and $v_5$. Since $n \ge 6$, we have $n' = n - 2 \ge 4$. We note that $H$ is also a mop. Then by the induction hypothesis, let $[D'_1,D'_2]$ be an LD-partition of $H$. In the following, we consider four cases and in each case, we prove that $G$ admits an LD-partition.

\medskip
\noindent\textbf{Case 2.1:} $v_1,v_2\in D'_1$ and $v_3\in D'_2$.

Let $D_1=D'_1\cup\{v_5\}$ and $D_2=D'_2\cup\{v_4\}$. Now we show that each $D_i$ is an LD-set of $G$ for $i\in\{1,2\}$. Note that $N_G(v_3)\cap\{v_1,v_2,v_5\}=\{v_1,v_2\}$ and $N_G(v_4)\cap\{v_1,v_2,v_5\}=\{v_1,v_2,v_5\}$. Moreover, $\{v_1,v_2,v_5\}\subset D_1$. Therefore, $D_1$ is an LD-set of $G$. Since $D'_2$ is an LD-set of $H$, we have $N_H(v_1)\cap D'_2\neq N_H(v_2)\cap D'_2$. Hence $N_G(v_1)\cap D_2\neq N_G(v_2)\cap D_2$. Note that $\{v_3,v_4\}\subseteq N_G(v_1)\cap D_2$, $\{v_3,v_4\}\subseteq N_G(v_2)\cap D_2$, and $N_G(v_5)\cap D_2=\{v_4\}$. Therefore, $D_2$ is an LD-set of $G$.

\medskip
\noindent\textbf{Case 2.2:} $v_1,v_3\in D'_1$ and $v_2\in D'_2$.

\medskip
Let $D_1=D'_1\cup\{v_4\}$ and $D_2=D'_2\cup\{v_5\}$. Now we show that each $D_i$ is an LD-set of $G$ for $i\in\{1,2\}$. Note that $N_G(v_2)\cap\{v_1,v_3,v_4\}=\{v_1,v_3,v_4\}$ and $N_G(v_5)\cap\{v_1,v_3,v_4\}=\{v_4\}$. Moreover, $\{v_1,v_3,v_4\}\subset D_1$. Therefore, $D_1$ is an LD-set of $G$. Since $D'_2$ is an LD-set of $H$, we have $N_H(v_1)\cap D'_2\neq N_H(v_3)\cap D'_2$. Hence $N_G(v_1)\cap D_2\neq N_G(v_3)\cap D_2$. Note that $v_5\notin N_G(v_1)\cap D_2$, $v_5\notin N_G(v_2)\cap D_2$, and $N_G(v_4)\cap D_2=\{v_2, v_5\}$. Therefore, $D_2$ is an LD-set of $G$.

\medskip
\noindent\textbf{Case 2.3:} $v_2,v_3\in D'_1$ and $v_1\in D'_2$. 

\medskip
Let $D_1=D'_1\cup\{v_5\}$ and $D_2=D'_2\cup\{v_4\}$.  We can make similar arguments as we did in Case 2.2.

\medskip
\noindent\textbf{Case 2.4:} $v_1,v_2,v_3\in D'_1$ and $v_1,v_2,v_3\notin D'_2$.

\medskip 
Let $D_1=D'_1\cup\{v_4\}$ and $D_2=D'_2\cup\{v_5\}$. Clearly, $D_1$ is an LD-set of $G$. Now we show that $D_2$ is an LD-set of $G$. Since $D'_2$ is an LD-set of $H$, $N_H(v_1)\cap D'_2$, $ N_H(v_2)\cap D'_2$, and $N_H(v_3)\cap D'_2$ are distinct and nonempty sets. Hence $N_G(v_1)\cap D_2$, $N_G(v_2)\cap D_2$, and $N_G(v_3)\cap D_2$ are distinct and nonempty sets. Note that $N_G(v_4)\cap D_2=\{v_5\}$. Moreover, we have $N_G(v_i)\cap D_2\neq N_G(v_4)\cap D_2$ for all $i\in\{1,2,3\}$ since $D'_2$ is an LD-set of $H$. Therefore, $D_2$ is an LD-set of $G$.
\end{proof}

\section{Split graphs and co-bipartite graphs}\label{sec:spco}

In this section, we show that every isolate-free and twin-free split graph and co-bipartite graph admit LD-partitions. In \cite{Foucaud16-2}, it has been proved that if $G$ is a twin-free and isolate-free split graph or co-bipartite graph, then  $\gamma_{LD}(G)\leq \frac{n}{2}$. Our proofs are an extension of the ones from~\cite{Foucaud16-2} used to prove Conjecture~\ref{conj} for split and co-bipartite graphs, but the key arguments are similar. 
 
\begin{theorem}
	If $G$ is an isolate-free and  twin-free split graph, then $G$ admits an LD-partition.
\end{theorem}
\begin{proof}
	Assume that $G$ is an  isolate-free and twin-free split graph. Let $X$ be  a clique of $G$  and $Y$ be an independent set of $G$ such that $X\cap Y=\emptyset$ and $X\cup Y=V(G)$. Note that every vertex in $Y$ has at least one neighbor in $X$. Moreover, for every $a,b\in X$, we have $N_G(a)\cap Y\neq N_G(b)\cap Y$ since $G$ is twin-free. Further, for every $c,d\in Y$, we have $N_G(c)\cap X\neq N_G(d)\cap X$.

    We may assume that the set $S = \{a \in X \mid N_G(a) \cap Y = \emptyset\}$ is empty (otherwise, $S$ has exactly one vertex, say $S = \{u\}$, and letting $X' = X \setminus \{u\}$ and $Y' = Y \cup \{u\}$, we obtain that the set $S'=\{a \in X' \mid N_G(a) \cap Y' = \emptyset\}$ is empty). Under this assumption, $[X, Y]$ is already an LD-partition of $G$. Thus $G$ admits an LD-partition.
\end{proof}

\begin{theorem}
	If $G$ is an isolate-free and  twin-free co-bipartite graph, then $G$ admits an LD-partition.
\end{theorem}
\begin{proof}
	Assume that $G$ is an isolate-free and twin-free co-bipartite graph. Let $X$ and $Y$ be the two cliques of $G$ such that $X\cap Y=\emptyset$ and $X\cup Y=V(G)$. Let $S_1=\{a\in X| N_G(a)\cap Y=\emptyset\}$ and $S_2=\{c\in Y| N_G(c)\cap X=\emptyset\}$. Since $G$ is twin-free, we have $N_G(a)\cap Y\neq N_G(b)\cap Y$ and $N_G(c)\cap X\neq N_G(d)\cap X$ for every $a,b\in X$ and $c,d\in Y$. So $|S_1|\le 1$ and $|S_2|\le 1$.
	
	If $S_1=\emptyset$ and $S_2=\emptyset$, then $X$ and $Y$ are LD-sets of $G$ and so $G$ admits an LD-partition. Suppose now that at least one of the sets from $S_1$ and $S_2$ is nonempty. First suppose that $|S_1|=1$ and $|S_2|=1$. Let $X'=(X\setminus S_1)\cup S_2$ and $Y'=(Y\setminus S_2)\cup S_1$. Now we show that $X'$ and $Y'$ are LD-sets of $G$. It is easy to observe that $X'$ is a dominating set of $G$. Since $N_G(c)\cap X\neq N_G(d)\cap X$ for every $c,d\in Y$, $N_G(c)\cap X'\neq N_G(d)\cap X'$ for every $c,d\in Y\setminus\{S_2\}$. Since $|S_2|=1$, we have $N_G(x)\cap X'\neq N_G(y)\cap X'$ for $x\in S_1$ and every $y\in Y\setminus\{S_2\}$. Hence $X'$ is an LD-set of $G$. Similarly, $Y'$ is also an LD-set of $G$ and so $G$ admits an LD-partition. Hence we assume that either $S_1=\emptyset$ or $S_2=\emptyset$. Without loss of generality, assume that $|S_1|=1$ and $S_2=\emptyset$. Let $S_1=\{x\}$. If there is no vertex $y\in Y$ such that $N_G(y)\cap X\setminus\{x\}= X\setminus\{x\}$, then let $X'=X\setminus\{x\}$ and $Y'=Y\cup \{x\}$. Then $X'$ and $Y'$ are LD-sets of $G$ and so $G$ admits an LD-partition. Hence there exists a vertex $y\in Y$ such that $N_G(y)\cap X\setminus\{x\}= X\setminus\{x\}$. Note that there is no vertex $y'\in Y$ other than $y$ such that $N_G(y')\cap X\setminus\{x\}= X\setminus\{x\}$; otherwise, $G$ contains twins. Let $X'=(X\setminus\{x\})\cup \{y\}$ and $Y'=(Y\setminus\{y\})\cup \{x\}$. Clearly, $X'$ is an LD-set of $G$. Now we show that $Y'$ is an LD-set of $G$. Note that $N_G(y)\cap Y'=Y'\setminus\{x\}$ and $N_G(v)\cap Y'\neq \emptyset$ for all $v\in X\setminus\{x\}$. Moreover, every vertex in  $X\setminus\{x\}$ has distinct neighborhood in $Y'$ since $G$ is twin-free and $x\in Y'$. Hence $Y'$ is also an LD-set of $G$. Thus $G$ admits an LD-partition.
\end{proof}

\section{Conclusion}\label{sec:con}

Conjecture~\ref{conj} is already known to hold for several important graph classes, including bipartite graphs, split graphs, co-bipartite graphs, line graphs, maximal outerplanar graphs, subcubic graphs, and block graphs. In this work, we addressed Question~\ref{motivation} and provided a positive answer for well-structured graph classes, namely distance-hereditary graphs, maximal outerplanar graphs, split graphs, and co-bipartite graphs.

This work contributes to a deeper structural understanding of locating-dominating sets and their properties within well-defined graph classes. As a natural direction for future research, it would be worthwhile to investigate whether Conjecture~\ref{conj} holds for the class of chordal graphs.

\section*{Acknowledgements}

Research of Florent Foucaud and Kaustav Paul was partially funded by the ANR project GRALMECO (ANR-21-CE48-0004), the French government IDEX-ISITE initiative 16-IDEX-0001 (CAP 20-25), the International Research Center ``Innovation Transportation and Production Systems'' of the I-SITE CAP 20-25, and the CNRS IRL ReLaX.

\section*{Declarations}

\noindent{\bf Conflict of interest} The authors do not have any financial or non financial interests that are directly or indirectly related to the work submitted for publication.

\noindent{\bf Data availability}
No data was used for the research described in this paper.

 \end{document}